

\magnification=\magstep1
\voffset=1.00truein
\settabs 18 \columns
\hoffset=1.00truein
\baselineskip=17 pt

\def\s{\smallskip}

\def\b{\bigskip}
\def\bb{\bigskip\bigskip}
\def\bbb{\bigskip\bigskip\bigskip}

\def\sqr#1#2{{\vcenter{\vbox{\hrule height.#2pt
 \hbox{\vrule width.#2pt height#1pt \kern#1pt
 \vrule width.#2pt} \hrule height.#2pt}}}}

\def\perp{\hbox{${\kern+.25em{\bigcirc}
\kern-.85em\bot\kern+.85em\kern-.25em}$}}
\def\lsim{\;raise0.3ex\hbox{$>$\kern-0.75em\raise-1.1ex\hbox{$\sim$}}\;}
\def\gsim{\;raise0.3ex\hbox{$>$\kerm-0.75em\raise-1.1ex\hbox{$\sim$}}\;}
\def\no{\noindent}

\def\ce{\centerline}
\def\ve{\vfill\eject}
\def\rdots{\mathinner{\mkern1mu\raise1pt\vbox{\kern7pt\hbox{.}}\mkern2mu
 \raise4pt\hbox{.}\mkern2mu\raise7pt\hbox{.}\mkern1mu}}

\def\e e{$e^+ e^-$ }


\rightline {UCLA/95/TEP/1}
\rightline {January 1995}
\bbb
\ce {\bf $q$-DEFORMATION OF THE LORENTZ GROUP}
\bb
\ce {Robert J. Finkelstein}
\b
\ce {\it Department of Physics}
\ce {\it University of California, Los Angeles, CA 90024-1547}
\bbb

\no {\bf ABSTRACT.}  We describe a $q$-deformation of the Lorentz
group in terms of a $q$-\break
deformation of the van der Waerden spinor algebra.
\ve

\line {\bf 1. Introduction.  \hfil}
\s

Weak deformations which test the stability of quantum mechanics
are of interest as a matter of principle.  One deformation which
has been the subject of considerable study results from the
replacement of the commutators or anticommutators of dynamically
conjugate operators by their $q$-commutators.  As far as these
studies have shown, there is no obstruction to the formulation of
a ``$q$-quantum mechanics" for finite systems.$^1$  In the field
theoretic context, however, this is not a well defined procedure
and may lead to violation of the Poincar\'e group.

In an alternative approach, one may begin with an explicit
deformation of this group.  We shall here examine the $q$-deformation
of the Lorentz group by $q$-deforming its two-dimensional
representation.  Although the $q$-Lorentz group has already been
studied by several authors,$^2$ the present treatment may be of
interest as a simple modification of the van der Waerden calculus.

\vskip.5cm

\line {\bf 2. Spinor Algebra for Lorentz Group.  \hfil}
\s

We shall first summarize the familiar spinor algebra for the Lorentz
group.  Let

$$\sigma^k=(\sigma^o,\vec\sigma) \qquad  \sigma^o=
  \left(\matrix{1 & 0\cr 0 & 1\cr}\right) \eqno(2.1)$$

\no and

$$X=x_k\sigma^k=\left(\matrix{x_o+x_3 & x_1-ix_2\cr
  x_1+ix_2 & x_o-x_3\cr}\right)~. \eqno(2.2)$$

\no Then

$$\hbox{det} X=x_o^2-\vec x^2=\eta^{k\ell}x_kx_\ell \eqno(2.3)$$
$$\eta^{k\ell}=\left(\matrix{1 & \hfil & \hfil & \hfil\cr
  \hfil & -1 & \hfil & \hfil\cr \hfil & \hfil & -1 & \hfil\cr
  \hfil & \hfil & \hfil & -1\cr}\right)~. $$

\no $X$ is also Hermitian:

$$X=X^\dagger~. \eqno(2.4)$$

\no Let

$$X^\prime=LXL^\dagger~. \eqno(2.5)$$

\no Then

$$(X^\prime)^\dagger=X^\prime~. \eqno(2.6)$$

\no But we shall now require

$$\hbox{det} L=1 \eqno(2.7)$$

\no then

$$\hbox{det} X^\prime=\hbox{det}X \eqno(2.8a)$$

\no or

$$(x_o^2-\vec x^2)^\prime=x_o^2-\vec x^2~. \eqno(2.8b)$$

\no $L$ has four complex matrix elements related by $\hbox{det}L=1$.
These are the six independent real parameters needed to describe
a Lorentz transformation and $L$ is a two-dimensional representation
of the Lorentz group.

Let $\xi^A$ be a two-rowed basis for the two-dimensional representation,
$L$,

$$\xi^A \rightarrow L^A_{~B}\xi^B \eqno(2.9)$$

\no and $\eta$ a basis for the conjugate representation

$$\eta^{\dot A}\rightarrow \bar L^{\dot A}_{~\dot B}
  \eta^{\dot B} \eqno(2.10)$$

\no where the conjugate representation is distinguished by a dotted
index.  Likewise let $\xi$ be a two-rowed basis for $L^{-1}$

$$\xi_A=(L^{-1})^{~B}_A\xi_B~. \eqno(2.11)$$

\no Again there is a conjugate representation:

$$\chi_{\dot A}\rightarrow (\bar L^{-1})_{\dot A}^{~\dot B}
  \chi_{\dot B}~. \eqno(2.12)$$

\no Then upper and lower indices behave in the usual way as
contravariant and covariant, so that expressions like
$\xi^A\xi_A$ are invariant.

The two-dimensional Levi-Civita symbol, $\epsilon_{AB}$, may be
used to define the determinant of $L$:

$$\epsilon_{AB}~ \hbox{det}~ L=L^C_{~A} L^D_{~B} \epsilon_{CD}~. \eqno(2.13)$$

\no Since $\hbox{det}L=1$,

$$\epsilon_{AB}=L^C_{~A} L^D_{~B} \epsilon_{CD}~. \eqno(2.14)$$

\no Therefore $\epsilon_{AB}$ is a covariant two-spinor that is
taken into itself by spin transformations while

$$\epsilon_{AB}\xi^A\eta^B \eqno(2.15)$$

\no is invariant.

Let us define $\epsilon^{AB}$ by

$$\epsilon^{AB}=\epsilon_{AB}~. \eqno(2.16)$$

\no Then

$$\epsilon^{AC}\epsilon_{BC}=\delta^A_B~. \eqno(2.17)$$

\no Also

$$\epsilon^{AB}=(L^{-1})^A_{~C}~(L^{-1})^B_{~D}~\epsilon^{CD}~. \eqno(2.18)$$

\no The spin tensors $\epsilon_{AB}$ and $\epsilon^{AB}$ provide ways
of lowering and raising indices.  Thus

$$\xi_B=\xi^A \epsilon_{AB} \eqno(2.19)$$

$$\xi^A=\epsilon^{AB} \xi_B \eqno(2.20)$$

\no and

$$\xi_A \eta^A=\xi^B\epsilon_{BA}\eta^A=
  -\eta^A \epsilon_{AB} \xi^B=
  -\eta_B \xi^B \eqno(2.21)$$

\no if $\xi$ and $\eta$ commute.  Therefore $\epsilon_{AB}$ serves
as a two-dimensional metric.

One next introduces the matrices contravariant to
$(\sigma^m_{A\dot x})$ with respect to $\epsilon$

$$(\bar\sigma^m)^{\dot XA}=
  \epsilon^{\dot X\dot Y} \epsilon^{AB}
  (\sigma^m)_{B\dot Y}~. \eqno(2.22)$$

\no Then if $(\sigma^m)_{B\dot Y}$ is the original set $(1,\vec\sigma)$,
defined in (2.1),

$$(\bar\sigma^m)^{\dot XA}=(1,-\vec\sigma)~. \eqno(2.23)$$

\no The following relations are also useful

$$\eqalignno{\bar\sigma^m \sigma^n+\sigma^n\bar\sigma^m &=
  2\eta^{nm} & (2.24) \cr
  \hbox{Tr}~\bar\sigma^m\sigma^n &= 2\eta^{nm} & (2.25) \cr}$$

\no or

$$\eqalignno{(\bar\sigma^m)^{\dot XA}
  (\sigma_n)_{A\dot X}&=2\delta^m_n & (2.26) \cr
  \noalign{\hbox{and}}
  (\sigma_n)_{A\dot X}(\sigma^n)^{B\dot Y} &= 2\delta^B_A
  \delta^{\dot Y}_{\dot X} & (2.27)\cr}$$

With the aid of the $\sigma$-matrices one may pass between the
four-dimensional and spin representation

$$T^{A\dot X~B\dot Y\ldots}_{\qquad C\dot W\ldots}=
  \sigma_a^{A\dot X} \sigma_b^{B\dot Y}
  \sigma_{C\dot W}^c T^{ab\ldots}_{\quad c\ldots} \eqno(2.28)$$

\no and

$$T^{ab\ldots}_{\quad c\ldots}=\sigma^a_{A\dot X}
  \sigma^b_{B\dot Y} \sigma^{C\dot W}_c
  T^{\dot A\dot X~B\dot Y\ldots}_{\qquad C\dot W\ldots}~. \eqno(2.29)$$

Any finite irreducible representation of $L$ is equivalent to some
spin representation which is separately symmetric in all dotted and
undotted indices.  $D(k,\ell)$ is the usual notation to describe
$2k$ dotted and $2\ell$ undotted indices.
\vskip.5cm

\line{\bf 3. Spin Representation of a Vector. \hfil}
\s

In the notation just introduced, the vector $x^a$ has the spin
representation $X^{A\dot X}$ where

$$X^{A\dot X} = \sigma_a^{A\dot X} x^a \eqno(3.1)$$

$$2x^a=\sigma^a_{A\dot X} X^{A\dot X}~. \eqno(3.2)$$

\no If

$$X^{A\dot X}=\psi^A\bar\psi^{\dot X} \eqno(3.3)$$

\no then

$$2x^a=\psi\sigma^a\bar\psi~. \eqno(3.4)$$

\no The transformation of $X^{A\dot X}$ is given by

$$\eqalignno{(X^{A\dot X})^\prime  &= \psi^{\prime A}
  \bar\psi^{\prime\dot X} & (3.5) \cr
  &= (L^A_{~B}\psi^B)(\bar L^{\dot X}_{~\dot Y} \bar\psi^{\dot Y}) \cr
  &= L^A_{~B}(\psi^B\bar\psi^{\dot Y})(L^\dagger)_{\dot Y}^{~\dot X} \cr
  &= L^A_{~B} X^{B\dot Y}(L^\dagger)_{\dot Y}^{~\dot X} \cr
  \noalign{\hbox{or}}
  X^\prime &= LXL^\dagger & (3.6a) \cr
  \noalign{\hbox{where}}
  X &= (X^{A\dot X}) = x_a\sigma^a & (3.7)\cr} $$

\no so that we recover (2.5).

{}From (3.6) one may obtain the vector representation of $L$ in the
following familiar way.  We have

$$\sum x^\prime_a \sigma^a=\sum(L\sigma^b L^\dagger)x_b \eqno(3.8)$$

$$x^\prime_b={1\over 2} \sum(\hbox{Tr}~\bar\sigma_b L\sigma^c L^\dagger)x_c
  \eqno (3.9)$$

\no or the vector representation of $L$ is

$$L_b^{~c}={1\over 2} \hbox{Tr}~ \bar\sigma_b L \sigma^c L^\dagger
  \eqno(3.10)$$

\no in terms of the spin representation of $L$.
\vskip.5cm

\line {\bf 4. The $q$-Deformation of the Lorentz Group. \hfil}
\s

We base our work on the following $q$-covariant
Levi-Civita tensor:

$$\epsilon_q=\left(\matrix{0 & q^{-{1\over 2}} \cr
  -q^{{1\over 2}} & 0 \cr} \right)   \quad
  \eqalign{& \epsilon^2_q=-1 \cr
  & \epsilon_q\epsilon^t_q=
  \left(\matrix{q^{-1} & 0\cr 0 & q\cr}\right)\cr} \eqno(4.1)$$

\no Repeating (2.13) in terms of $\epsilon_q$ we define $T$, the $q$-
deformed spin representation of $L$ as follows:

$$\eqalignno{\epsilon_q\Delta_q &= T^t \epsilon_q T & (4.2) \cr
  &= T \epsilon_q T^t & (4.3) \cr}$$

\no Here $t$ means transpose, and $\Delta_q$ is by definition
the $q$-determinant.

We know that these relations define the $GL_q(2)$ group.
That is, if

$$T=\left(\matrix{a & b\cr c & d\cr}\right) \eqno(4.4)$$

\no then (4.2) and (4.3) may be simultaneously satisfied only if
the following relations hold:

$$\eqalign{ab &= qba\cr  ac &= qca\cr  bd &= qdb\cr  cd &= qdc \cr}
  \qquad
  \eqalign{bc &= cb\cr  \hfil \cr  \Delta_q &= ad-qbc \cr
  \hfil \cr} \eqno(4.5)$$

\no and

$$(\Delta_q,T_{k\ell})=0~. \eqno(4.6)$$

\no $\Delta_q$ is the $q$-determinant and (4.2) may be regarded
as a definition of this determinant since:

$$\eqalignno{(\epsilon_q)_{AB}~\hbox{det}_q~T &=
  T^C_{~A}T^D_{~B}(\epsilon_q)_{CD} & (4.7) \cr
  &= T_A^{~C} T^{~D}_B (\epsilon_q)_{CD} & (4.8) \cr}$$

\no as in (2.13).

The only restriction on the matrix $L$ (first introduced in
(2.3)) is (2.7).  We shall now impose the corresponding
condition on $T$.  Then $\hbox{det}_q~T=1$ and

$$\eqalignno{\epsilon_q &= T\epsilon_q T^t & (4.9)\cr
  &= T^t \epsilon_q T & (4.10)\cr} $$

\no Note also

$$\hbox{det}_q T^t=\hbox{det}_qT=1~. \eqno(4.11)$$

We may now try to take over Eq. (2.5) with $L$ replaced by $T$

$$X^\prime=TXT^\dagger \quad \hbox{det}_qT=1~. \eqno(4.12)$$

The previous arguments no longer go through, however, since the
matrix elements of $T$ are non-commuting.  In the limit $q=1$, $T$
approaches $L$, and (4.12) then describes a Lorentz transformation.
However, (4.11) has the same meaning as (2.5) only in the $q=1$ limit
when the matrix elements $(a,b,c,d)$ all commute and lie in the
complex plane.  In general (4.12) will have no classical meaning.
Nevertheless a corresponding spinor algebra may again be constructed.
\vskip.5cm

\line {\bf 5. $q$-Spinors. \hfil}
\s

Except when explicitly otherwise indicated, let us now understand
by $\epsilon$ the $\epsilon_q$ matrix.  One may next define a
contravariant $\epsilon$ by

$$\epsilon^{12}=q^{{1\over 2}}~, \quad
  \epsilon^{21}=-q^{-{1\over 2}}~, \quad
  \epsilon^{kk}=0 \eqno(5.1)$$

\no so that $\epsilon^{AB}(q)=\epsilon_{AB}(q^{-1})$ or

$$\epsilon^{AC}\epsilon_{BC} = \delta^A_B~. \eqno(5.2)$$

\no Then $\epsilon_{AB}$ is the metric in spin space and
may be used to raise and lower indices.  Thus if $\xi_A$ is
a covariant spinor, the corresponding contravariant spinor is

$$\xi^B=\epsilon^{BA}\xi_A \eqno(5.3)$$

\no and

$$\xi_A=\xi^C\epsilon_{CA}~. \eqno(5.4)$$

\no By (4.7) $\epsilon_{AB}$ is an invariant tensor since

$$\epsilon_{AB}=T^{~C}_A T^{~D}_B \epsilon_{CD}~. \eqno(5.5)$$

\no Multiply (5.5) by $\xi^A \chi^B$.  Then

$$\xi^A\epsilon_{AB}\chi^B = \xi^{C^\prime}\epsilon_{CD}
  \chi^{D^\prime} \eqno(5.6)$$

\no where

$$\eqalignno{\xi^{C^\prime} &= \xi^AT^{~C}_A & (5.7a)\cr
  \chi^{D^\prime} &= \chi^B T^{~D}_B~. & (5.7b)\cr}$$

\no The invariant quadratic forms (5.6) may also be written as
follows:

$$\xi^A\epsilon_{AB}\chi^B=\xi_B\chi^B=\xi^A\tilde\chi_A \eqno(5.8)$$

\no where

$$\eqalignno{\xi_B &= \xi^A\epsilon_{AB} & (5.9a)\cr
   \tilde\chi_A &= \epsilon_{AB}\chi^B=\chi^B\epsilon^t_{BA} & (5.10a)\cr}$$

\no By (5.2) the inverses of (5.9a) and (5.10a) are

$$\eqalignno{\xi^C &= \epsilon^{CB}\xi_B & (5.9b)\cr
  \chi^C &= \tilde\chi_A\epsilon^{AC} & (5.10b)\cr}$$

\no Given (5.7) and (5.9a) one finds

$$\eqalignno{\xi^\prime_B &= \xi^CT^{~A}_C\epsilon_{AB} & (5.11)\cr
  \noalign{\hbox{By (5.9b)}}
  \xi^\prime_B &= \epsilon^{CD}\xi_D T^{~A}_C\epsilon_{AB} & (5.12)\cr}$$

\no Denote the matrix $\parallel \epsilon^{CD}\parallel$ by
$\hat\epsilon$ to distinguish it from $\parallel \epsilon_{CD}\parallel
=\epsilon$.  Then

$$\xi^\prime_B=\xi_D(\hat\epsilon^t T\epsilon)^D_{~B}~. \eqno(5.13)$$

But

$$\hat\epsilon(q)=\epsilon(q^{-1})=-\epsilon^t(q) \eqno(5.14)$$

\no and (5.13) becomes

$$\xi_B^\prime=-\xi_D(\epsilon T\epsilon)^D_{~B}~. \eqno(5.15)$$

\no By (4.9) or (4.10)

$$\epsilon T\epsilon =-(T^t)^{-1}~. \eqno(5.16)$$

\no Therefore

$$\xi^\prime_B=\xi_D\bigl((T^t)^{-1}\bigr)^D_{~B}~. \eqno(5.17)$$

\no Note that $(T^t)^{-1}\not= (T^{-1})^T$ here.  One checks that

$$\xi^\prime_B\chi^{B^\prime}=\xi_B\chi^B~. \eqno(5.18)$$

\no This invariant may be written as

$$\eqalignno{\xi_A\chi^A &= \epsilon_{BA}\xi^B \chi^A & (5.19)\cr
  &= q^{-{1\over 2}}(\xi^1\chi^2-q\xi^2\chi^1)~. & (5.20)\cr}$$

\no Therefore the invariant equation

$$\xi_A\chi^A=0 \eqno(5.21)$$

\no implies the invariance of the commutation rules

$$\xi^1\chi^2=q\xi^2\chi^1~. \eqno(5.22)$$

Just as in the Lorentz case one must also make use of the conjugate
representations.  Corresponding to (5.7) and (5.17) one has

$$\eqalignno{\bar\xi^{\dot C^\prime} &= \bar\xi^{\dot A}
  \bar T_{\dot A}^{~\dot C} & (5.23)\cr
  \noalign{\hbox{and}}
  \bar\xi^\prime_{\dot B} &= \bar\xi_{\dot D}
  \bigl((\bar T^t)^{-1}\bigr)^{\dot D}_{~\dot B} & (5.24)\cr}$$

\no where $\bar\xi^{\dot A}$ is the basis for $\bar T$, the
conjugate representation, and the dot again indicates the
conjugate representation.  The conjugate is now in the $SU_q(2)$
algebra, not in the complex plane.

An alternative procedure begins with (4.7) instead of (4.8).  Then

$$\eqalign{&\epsilon_{AB}=T^C_{~A}\epsilon_{CD}T^D_{~B} \cr
  & \xi^A\epsilon_{AB}\chi^B=(\xi^C\epsilon_{CD}\chi^D)^\prime \cr}$$

\no where

$$\eqalign{\xi^{C^\prime} &= T^C_{~A}\xi^A \cr
  \chi^{D^\prime} &= T^D_{~B}\chi^B \cr}$$

\no instead of (5.7).

\vskip.5cm

\line {\bf 6. The $\sigma_q$ Matrices. \hfil}
\s

Let us set

$$(\sigma^m_q)_{B\dot Y} = (1,\vec\sigma) \eqno(6.1)$$

\no just as for the Lorentz case.

We now introduce the matrices contravariant to $(\sigma^m_q)_{B\dot Y}$
with respect to the metric $\epsilon_q$:

$$(\bar\sigma^m_q)^{\dot XA}=
  \epsilon_q^{\dot X\dot Y} \epsilon_q^{AB}
  (\sigma_q^m)_{B\dot Y}~. \eqno(6.2)$$

\no Then

$$(\bar\sigma^m_q)^{\dot XA}=
  \left(\matrix{q & 0\cr  0 & q^{-1}\cr}\right)~
  \left(\matrix{0 & -1\cr -1 & 0\cr}\right)~
  \left(\matrix{0 & i\cr -i & 0\cr}\right)~
  \left(\matrix{-q & 0 \cr 0 & q^{-1}\cr}\right) \eqno(6.3)$$

\no which satisfy the following relations

$$\eqalignno{(\bar\sigma_q^m)^{\dot XA}
  (\sigma_q^n)_{A\dot X} &= 2\eta^{mn} & (6.4) \cr
  (\sigma_q^n)_{A\dot X}(\bar\sigma_{qn})^{\dot YB} &=
  2\delta^{\dot Y}_{\dot X} \delta^B_A & (6.5)\cr}$$

\no where

$$\eta^{nm}=
  \left(\matrix{{1\over 2}(q+q^{-1}) & 0 & 0 & {1\over 2}(q-q^{-1})\cr
  0 & -1 & 0 & 0 \cr
  0 & 0 & -1 & 0 \cr
  {1\over 2}(q-q^{-1}) & 0 & 0 & -{1\over 2}(q+q^{-1})\cr}\right)~.\eqno(6.6)$$
\vskip.5cm

\line {\bf 7. The $q$-Spin Representation of Vectors. \hfil}
\s

Define the $q$-vector by

$$x_a=(\sigma_q)_a^{A\dot X} X_{A\dot X} \eqno(7.1)$$

\no where $X_{A\dot X}$ is the bispinor:

$$X_{A\dot X}=\psi_A \bar\psi_{\dot X}~. \eqno(7.2)$$

\no Then

$$X_a=\psi(\sigma_q)_a\bar\psi~. \eqno(7.3)$$

\no The transformation of $X_{A\dot X}$ reads as follows:

$$\eqalign{X^\prime_{A\dot X} &= (T_A^{~B} \psi_B)
  (\bar T_{\dot X}^{~\dot Y}\bar\psi_{\dot Y}) \cr
  &= T_A^{~B}(\psi_B\bar\psi_{\dot Y})
  (T^\dagger)^{\dot Y}_{~\dot X} \cr} \eqno(7.4)$$

\no where

$$\bar T^{~\dot Y}_{\dot X} = (\bar T^t)^{\dot Y}_{~\dot X}=
  (T^\dagger)^{\dot Y}_{~\dot X}~. \eqno(7.5)$$

\no Then (7.4) may be written

$$X^\prime=TXT^\dagger \eqno(7.6)$$

\no where $T^\dagger$ means the Hermitian adjoint matrix.  This equation is
formally the same as (3.6) with $L$ replaced by $T$.  In addition

$$\hbox{det}_qT=\hbox{det}~L=1 \eqno(7.7)$$

\no but

$$\hbox{det}_q X^\prime \not= \hbox{det}_qX \eqno(7.8)$$

\no unless $q=1$, since

$$\hbox{det}_q AB \not= \hbox{det}_qA~\hbox{det}_qB \eqno(7.9)$$

\no unless

$$(A_{ij},B_{k\ell})=0~.$$

\no Although $\hbox{det}_qX$ is not conserved, one may continue
by rewriting (7.6) as follows:

$$\sum x^{a^\prime}(\sigma_q)_a=T\sum x^a(\sigma_q)_a
  T^\dagger \eqno(7.10)$$

\no and solving to give the transformation equation for $x^a$

$$x^{a^\prime}=\sum T^a_{~b}x^b \eqno(7.11)$$

\no where

$$T^a_{~b} = {1\over 2} \hbox{Tr}~\bar\sigma^a_q
  T(\sigma_q)_b T^\dagger~. \eqno(7.12)$$
\vskip.5cm

\line {\bf 8. Invariants and Irreducible Representations.$^3$  \hfil}
\s

Although the usual Lorentz interval, represented by $(\hbox{det}X)^{1/2}$,
is no longer invariant, there are of course new invariants belonging
to the deformed group.  Just as the invariants of the Lorentz
group are the same as the invariants of the unimodular linear
group $SL(2)$, here they are the same as those of $SL_q(2)$.   The basic
invariant may be expressed in terms of either the covariant or
contravariant metric as follows:

$$\eqalignno{\xi^A\epsilon_{AB}\chi^B &=0  & (8.1a) \cr
  \xi_A \hat\epsilon^{AB}\chi_B &= 0 & (8.1b) \cr}$$

\no or

$$\eqalignno{\xi^1\chi^2-q\xi^2\chi^1 &=0 & (8.2a)\cr
  \xi_2\chi_1-q\xi_1\chi_2 &=0~. & (8.2b) \cr}$$

\no Take the special case $\xi=\chi$.  Then

$$\chi_2\chi_1=q\chi_1\chi_2 \eqno(8.3)$$

\no and

$$\chi_1\chi_2 \cdot \chi_2\chi_1=\chi_2\chi_1\cdot
  \chi_1\chi_2~. \eqno(8.4)$$

\no By the binomial theorem applied to (8.1)

$${(i\chi^t\epsilon \chi)^{2j}\over (2j)!} =
  \sum \tilde V(jm) V(jm) (-q)^m \eqno(8.5)$$

\no where

$$V(jm) =
  {\chi_2^{j+m} \chi_1^{j-m}\over
  [(j+m)!(j-m)!]^{1/2}}~\epsilon(q^{-1}|j+m) \eqno(8.6)$$

$$\tilde V(jm) =
  {\chi_1^{j+m}\chi_2^{j-m}\over
  [(j+m)!(j-m)!]^{1/2}}~
  \epsilon(q|j-m) \eqno(8.7)$$

\no and

$$\epsilon(q|n)=q^{n(n-1)/2}~.  \eqno(8.8)$$

\no Therefore

$$\eqalign{e^{i\chi^t\epsilon\chi}&= \sum
  \tilde V(jm)V(jm)(-q)^m \cr
  &=1 \cr}\eqno(8.9)$$

\no The invariant terms of (8.5) may be written

$$Q(j)=\tilde V(j) C^j(-q)V(j) \eqno(8.10)$$

\no where

$$C^j(q)=\left(\matrix{q^j & \hfil & \hfil \cr
  \hfil & \ddots & \hfil\cr
  \hfil & \hfil & q^{-j}\cr}\right)~. \eqno(8.11)$$

\no Here $\tilde V(jm)$ and $V(jm)$ are the $2j+1$ dimensional
vectors which transform under $T$ as

$$\eqalignno{\tilde V^\prime(jm) &= \sum^j_{-j}
  \tilde V(jm^\prime) \tilde D^j(m^\prime m) & (8.12)\cr
  V^\prime (jm) &= \sum^j_{-j} D^j(m,m^\prime)
  V(jm^\prime) & (8.13)\cr}$$

\no One finds

$$\eqalign{\tilde D^j(m^\prime,m) &= N^{jm}_{jm^\prime}
  {\epsilon(q|j-m)\over \epsilon(q|j-m^\prime)} q^{-j^2} \cr
  &\times \sum
  q^\sigma \biggl\langle\matrix{j+m\cr j-m^\prime-t\cr}
  \biggr\rangle_{q^2}~
  \biggl\langle\matrix{j-m\cr t\cr}\biggr\rangle_{q^2}~
  b^{j-m^\prime-t}a^{m+n^\prime+t}d^tc^{j-m-t}~. \cr} \eqno(8.14)$$

\no Here

$$\biggl\langle\matrix{n\cr s\cr}\biggr\rangle=
  {\langle n\rangle_q!\over
  \langle s\rangle_q!~\langle n-s\rangle_q!} \qquad
  \langle n_q\rangle = {q^n-1\over q-1} \eqno(8.15)$$

\no where $\langle n\rangle$ is the basic integer.

These $2j+1$ dimensional vectors are the higher dimensional spin
tensors that correspond to the higher rank Lorentz tensors, and
the invariants $Q(j)$ replace the Lorentz invariants such as the
interval $(\hbox{det}X)^{1/2}$.

To put (8.6) into correspondence with the bispinor of (7.2),
one sets $j=1/2$ in (8.6) to obtain $\psi^A$, and then makes
use of the conjugate representation in order to form
$\chi^A\bar\chi^{\dot X}$.

Although it is true that one recovers the Lorentz group when $q$
is set equal to unity, the results that hold when $q\not= 1$ are
totally different in nature from the $q=1$ formulas and can have
significance only in a quantized theory.
The parameter $q$ appearing in this paper is, of course, also
totally different from the $q$ appearing in $q$-commutators.

\ve

\line {\bf References. \hfil}
\s
\item{1.} E. Marcus and R. Finkelstein, UCLA/94/TEP/32 (1994),
J. Math. Phys. (in press).
\item{2.} R. Finkelstein, UCLA/94/TEP/46 (1994), Lett. Math. Phys.
(in press).
\item{3.} O. Ogievetsky, W. B. Shmidke, J. Wess, and B. Zumino,
Lett. Math. Phys. {\bf 23}, 233 (1991);
\item{} P. Truini and V. S. Varadarajan, Lett. Math. Phys.
{\bf 21}, 287 (1991); {\it ibid.} {\bf 26}, 53 (1992).
\item{4.} G. Lustig, Adv. Math. {\bf 20}, 237 (1988);
\item{} M. Rosso, Comm. Math. Phys. {\bf 117}, 581 (1988);
\item{} T. H. Koornwinder, NATO ASI Ser C{\bf 294}, 257 (1994);
\item{} R. Finkelstein, Lett. Math. Phys. {\bf 24}, 75 (1993).

\bye